\documentclass{amsart}
\usepackage{epsfig, amsfonts}

\newcommand{\N}{{\mathbb N}}
\newcommand{\Z}{{\mathbb Z}}
\newcommand{\e}{\varepsilon }
\newcommand{\G} {\Gamma }
\newcommand{\gr}{Grid _\delta }
\newcommand{\dc}{dist _{Cone}}
\newcommand{\dm}{dist _M}
\newcommand{\dg}{dist _G}
\newcommand{\diam}{diam \; }
\newcommand{\C}{{\rm Cone}_\omega (G)}
\newcommand{\card }{\sharp }
\newcommand{\Ca}{Cay (G)}

\newtheorem{thm}{Theorem}[section]
\newtheorem{cor}[thm]{Corollary}
\newtheorem{lem}[thm]{Lemma}
\newtheorem{prop}[thm]{Proposition}
\newtheorem{prob}[thm]{Problem}

\theoremstyle{definition}
\newtheorem{defn}[thm]{Definition}
\theoremstyle{remark}

\title{Fundamental groups of asymptotic cones.}

\author{A.Erschler, D.Osin}

\address{CNRS, University Lille 1, UFR de mathematiques, 59655, Villeneuve d'Ascq Cedex, France}
\email{erschler@math.univ-lille1.fr, erschler@pdmi.ras.ru}

\address{Department of Mathematics, 1326 Stevenson Center,
 Vanderbilt University, Nashville  TN 37240-0001, USA}
\email{denis.ossine@vanderbilt.edu}

\subjclass[2000]{Primary: 20F69; Secondary: 14F35, 20F06}
\keywords{asymptotic cone, isometric embedding at infinity, finitely
  generated group}
\thanks{This work was partially supported by
the Swiss National Science foundation. The second author was also
supported by the RFBR grant $\# $ 99-01-00894.}

\usepackage{amssymb, amscd, amsmath}

\begin{document}

\begin{abstract}
We show that for any metric space $M$ satisfying certain natural conditions, there is a finitely
generated group $G$, an ultrafilter $\omega $, and an isometric embedding $\iota $ of $M$ to the
asymptotic cone $\C $ such that the induced homomorphism $\iota ^ \ast :\pi_1(M)\to \pi_1(\C )$ is
injective. In particular, we prove that any countable group can be embedded into a fundamental group
of an asymptotic cone of a finitely generated group.
\end{abstract}

\maketitle

\section{Introduction}

To any metric space $X$ with a distance function $dist$, one can
associate a new metric space $Cone_\omega (X)$, the so--called
{\it asymptotic cone} of $X$, by taking the ultralimit of scaled
spaces $(X, \frac{1}{n}dist)$ with respect to an ultrafilter
$\omega $. Informally speaking $Cone_\omega (X)$ shows what $X$
looks like if the observer is placed at 'infinity' (see the next
section for precise definition). This notion appears in the proof
of Gromov's theorem about groups of polynomial growth (it is
defined in  \cite{wilk}, though in the polynomial growth case,
where for the convergence to the limit one does not need
ultralimits, the corresponding limit space is considered already
in \cite{grom}).

It is known that a group is hyperbolic if and only if any of its
asymptotic cones is an $R$-tree (\cite{groa}). (Misha Kapovich
pointed out to the first author that if one of the cones of a
 finitely presented group is an $R$-tree, then this group is hyperbolic).

 The application of asymptotic cones to the study of algebraic and geometric properties of
hyperbolic groups  appears in a different language in
\cite{Morgan} and  in \cite{rips}. For a recent progress see
\cite{Sela1,Sela2,Sela3} and references there. Asymptotic cones
can be used for proving rigidity theorems for symmetric spaces
\cite{KL}. For other results about asymptotic cones we refer to
\cite{brid,buri,drutu,jaos,papa,paulin, rile, thve,new}.

The case when $X$ is a finitely generated group $G$ endowed with a
word metric is of particular interest. In \cite{groa}, Gromov
pointed out a connection between homotopical properties of $\C $
and asymptotic invariants of $G$. Namely, if $\C $ is simply
connected with respect to all ultrafilters, then $G$ is finitely
presented and the Dehn function is polynomial.
 A partial converse result was obtained in
\cite{papa}. However almost nothing was known about algebraic
structure of the fundamental group $\pi _1(\C )$ in case $\pi
_1(\C )$ is nontrivial. In particular, no examples of non--free
finitely generated subgroups of such a fundamental groups were
known until now.  (It was observed in \cite{buri} that  the
asymptotic cones of Baumslag--Solitar groups
$$BS_{p,q}=\langle a,b\; | \; b^{-1}a^pb = a^q\rangle ,$$ where
$|p|\ne |q|$, contain non-free infinitely generated groups.)

The following question is stated in  \cite{groa}.

\begin{prob}\label{mp}
Which groups can appear as (finitely generated) subgroups in
fundamental groups of asymptotic cones of finitely generated
groups?
\end{prob}

In the present paper we answer this question  by proving the following theorem.

\begin{thm}\label{main}
Let $M$ be a metric space such that

\begin{enumerate}

\item[(M1)] $M$ is geodesic, that is for any two points $x,y\in M$
there is a path joining $x$ and $y$ of length $\dm (x,y)$;

\item[(M2)]  there is a sequence of compact subsets $M_1\subseteq
M_2\subseteq \ldots $ of $M$ such that
$M=\bigcup\limits_{i=1}^\infty M_i$.

\end{enumerate}

\noindent Then there exists a group $G$ generated by $2$ elements,
an ultrafilter $\omega $, and an isometric embedding $\iota $ of
$M$ into the asymptotic cone $\C $. If, in addition, $M$ is
uniformly locally simply--connected, i.e.,

\begin{enumerate}
\item[(M3)] there exists $\e >0$ such that
any loop in any ball of radius $\e$ is contractible,
\end{enumerate}
then the induced homomorphism of the fundamental groups $$\iota ^
\ast :\pi_1(M)\to \pi_1(\C )$$ is injective.
\end{thm}

Obviously any combinatorial complex $M$ with countable number of
cells in each dimension admits a (natural) metric which induces
the standard topology on $M$ and with respect to which $M$
satisfies (M1)--(M3). Since any countable group is a fundamental
group of a countable combinatorial 2-complex, we obtain

\begin{cor}
For any countable group $H$, there exists a finitely generated
group $G$ and an ultrafilter $\omega $ such that $\pi_1(\C )$
contains a subgroup isomorphic to $H$.
\end{cor}

The construction in the proof of Theorem \ref{main} is similar to
that in \cite{Ols} and can be intuitively understood as follows.
Given a metric space $M$ with metric $\dm $ satisfying (M1)--(M3),
we first approximate $M$ by finite $\e_i$--nets $Net _i$, where
$\e_i \to 0$ as $i\to \infty $. Then we choose a rapidly growing
sequence of natural numbers $\{ n_i\}$ and use a construction
similar to that from \cite{Ols} to produce embeddings $\alpha_{i}$
of $Net_{i}$ into a finitely generated group $G$ endowed with a
word metric $\dg $ such that $\alpha _{i}$, being considered as a
map from a metric space $(Net_{i}, \dm )$ to $(G, \frac{1}{n_i}\dg
)$, is a $(\lambda_i,c_i)$--quasi--isometry, where $\lambda _i\to
1$ and $c_i\to 0$ as $i\to \infty $. This gives an isometric
embedding $\iota :M\to \C $ for any ultrafilter $\omega $
satisfying $\omega \{ n_i\} =1$. Then condition (M3) allows to
show that $\iota $ induces an injective map on the fundamental
group $\pi _1(M)$.

The paper is organized as follows. In the next section we collect
main definitions and results which are used in what follows. Some
auxiliary results about words with small cancellation properties
are proved in Section 3. In Section 4 we construct the group $G$
and the injective map $\iota : M\to \C$ mentioned in Theorem 1 and
show that $\iota $ is an isometry.  In Section 5 we conclude the
proof of the main theorem by proving injectivity of  $\iota ^ \ast
:\pi_1(M)\to \pi_1(\C )$.

\section{Preliminaries}

\subsection{Asymptotic cones}

Recall that a {\it non-principal ultrafilter $\omega$ } is a
finitely additive non-zero measure on the set of all subsets of
${\Bbb N}$ such that each subset has measure either $0$ or $1$ and
all finite subsets have measure $0$. For any bounded function $h :
{\Bbb N}\to {\Bbb R}$ its limit $h(\omega)$ with respect to a
non-principal ultrafilter $\omega$ is uniquely defined by the
following condition: for every $\delta >0$
$$\omega(\{i\in \N \; :\;  |h(i)-h(\omega)|<\delta \})=1.$$

\begin{defn} Let $X$ be a metric space with a distance function $dist $.
We fix a basepoint $O\in X$ and consider the set of all sequences
$g:{\Bbb N}\to X$ such that
$$dist(O, g(i))\le \mbox{const}_g \cdot i$$ (here the constant
$const_g$ depends on $g$). To any pair of such sequences $g_1,g_2$
one may assign a function
$$h_{g_1,g_2}(i)=\frac{dist(g_1(i),g_2(i))}{i}.$$
We say that the sequences $g_1$, $g_2$ are equivalent if the limit
$h_{g_1,g_2}(\omega)=0$. The set $\C $ of all equivalence classes
of sequences with the distance
$$\dc (g_1,g_2)=h_{g_1,g_2}(\omega)$$ is called an {\it asymptotic
cone} of $X$ with respect to the non-principal ultrafilter
$\omega$. Clearly this space does not depend of the basepoint
chosen.
\end{defn}

If $G$ is a group generated by a finite set $S$, we can regard $G$
as a metric space assuming the distance between two elements
$a,b\in G$ to be equal to the length of the shortest word in the
alphabet $S^{\pm1}$ representing $a^{-1}b$. Such a distance
function is called the {\it word metric associated to $S$}. Given
an ultrafilter $\omega $, this leads to the asymptotic cone $\C $.
It is worth noting that $\C $ strongly depends on the choice of
the ultrafilter \cite{thve}.

\subsection{Cayley graphs and van Kampen diagrams}

Recall that a {\it Cayley graph} $\Ca $ of a group $G$ generated
by a set $S$ is an oriented labelled 1--complex with the vertex
set $V(\Ca )=G$ and the edge set $E(\Ca )=G\times S$. An edge
$e=(g,s)\in E(\Ca )$ goes from the vertex $g$ to the vertex $gs$
and has the label $Lab (e)=s$. As usual, we denote the origin and
the terminus of the edge $e$, i.e., the vertices $g$ and $gs$, by
$e_-$ and $e_+$ respectively. The word metric on $G$ associated to
$S$ can be extended to $\Ca $ by assuming the length of every edge
to be equal to one. Also, it is easy to see that a word $W$ in
$S^{\pm 1}$ represents $1$ in $G$ if and only if some (or,
equivalently, any) path $p$ in $\Ca $ labelled $W$ is a cycle.

A planar {\it map} $\Delta $ over a group presentation
\begin{equation}
G=\langle S\; | \; \mathcal P\rangle \label{ZP}
\end{equation}
is a finite oriented connected simply--connected 2--complex
endowed with a labelling function $Lab : E(\Delta )\to S^{\pm 1}$
(we use the same notation as for Cayley graphs) such that $Lab
(e^{-1})=(Lab (e))^{-1}$.

Given a combinatorial path $p=e_1e_2\ldots e_k$ in $\Delta $
(respectively in $\Ca $), where $e_1, e_2, \ldots , e_k\in
E(\Delta )$ (respectively $e_1, e_2, \ldots , e_k\in E(\Ca )$), we
denote by $Lab (p)$ its label. By definition, $Lab
(p)=Lab(e_1)Lab(e_2)\ldots Lab (e_k).$ We also denote by
$p_-=(e_1)_-$ and $p_+=(e_k)_+$ the origin and the terminus of $p$
respectively. A path $p$ is called {\it irreducible} if it
contains no subpaths of type $ee^{-1}$ for $e\in E(\Delta  )$
(respectively $e\in E(\Ca )$). The length $|p|$ of $p$ is, by
definition, the number $k$ of edges of $p$.

Given a cell $\Pi $ of $\Delta $, we denote by $\partial \Pi$ the
boundary of $\Pi $; similarly, $\partial \Delta $ denotes the
boundary of $\Delta $. The label of $\partial \Pi $ or $\partial
\Delta $ is defined up to a cyclic permutation. A map $\Delta $
over a presentation (\ref{ZP}) is called a {\it van Kampen
diagram} over (\ref{ZP}) if the following holds. For any cell $\Pi
$ of $\Delta $, the boundary label $Lab (\partial \Pi)$ is equal
to a cyclic permutation of a word $P^{\pm 1}$, where $P\in
\mathcal P$. Sometimes it is convenient to use the notion of
$0$--refinement in order to assume diagrams to be homeomorphic to
a disc. We do not explain here this notion and refer the
interested reader to \cite{Olsbook}.

The van Kampen lemma states that a word $W$ over the alphabet
$S^{\pm 1}$ represents the identity in the group given by
(\ref{ZP}) if and only if there exists a simply--connected planar
diagram $\Delta $ over (\ref{ZP}) such that $Lab (\partial \Delta
)$ coincides with $W$ \cite{LS}, \cite{Olsbook}. A van Kampen
diagram is called {\it minimal}, if it contains the minimal number
of cells among all diagrams with the same boundary labels.

\subsection{Approximations of metric spaces by graphs}

We recall that a subset $Z$ is called an $\delta $--net in a
metric space $Y$, if for all $y\in Y$ there exists $z\in Z$ such
that the distance between $y$ and $z$ is less than $\delta $. We
say that a subset $Z$ of a metric space $Y$ is a $(\delta _1,
\delta _2)$--{\it net} if $Z$ is a $\delta_1$--net and the
distance between any two points of $Z$ is greater than $\delta
_2$. The following lemma will be used in Section 4.

\begin{lem} \label{net1}
Suppose that $M$ is a metric space satisfying conditions (M1) and
(M2). There exists a sequence of finite subsets $Net_1\subseteq
Net_2\subseteq \ldots $ of $M$ such that for all $i\in \mathbb N$,
$Net_i$ is a $(2/i , 1/i)$--net in $M_i$.
\end{lem}

\begin{proof}
We proceed by induction on $i$. Suppose that $Net_{i-1}$ is an
$(2/(i-1) , 1/(i-1))$--net in $M_{i-1}$ if $i>1$, and
$Net_{i-1}=\emptyset $ if $i=1$. We consider an arbitrary finite
$1/i $--net $N$ in $M_i$ that contains $Net_{i-1}$ as a subset.
Let $\mathcal N$ denote the set of all subsets $L$ such that
$Net_{i-1}\subseteq L\subseteq N$ and for any two different
elements $x,y\in L$ we have $\dm (x,y)>1/i $. Note that the set
$\mathcal N$ is non--empty as it contains $Net _{i-1}$.

Consider a partial order on $\mathcal N$ which corresponds to
inclusion, i.e., for any $A,B\in \mathcal N$, $A\preceq B$ if and
only if $A\subseteq B$. Since $\mathcal N$ is finite, we can take
a maximal subset $B$ with respect to this order. Note that for any
$t \in N$, we have $\dm (t, B) \le 1/i $. Indeed, otherwise $B\cup
\{ t\} \in \mathcal N$ and thus $B$ is not maximal. Therefore, for
any $x\in M_i$, we have $\dm (x,B)\le 2/i $. Thus $B$ is an $(2/i,
1/i)$--net in $M_i$.
\end{proof}

\section{Words with small cancellations}

To prove the main result of our paper we will need an infinite set
of words satisfying a certain small cancellation conditions. We
begin with definitions.

Let $X$ be an alphabet and $F$ a free group with the basis $X$.
Throughout the following discussion we write $U\equiv V$ to
express the letter--by--letter equality of the words $U$ and $V$.
Given a word $W$ over the alphabet $X$, by a cyclic word $W$ we
mean the set of all cyclic shifts of $W$. Two cyclic words $W_1$
and $W_2$ are equal if and only there exist cyclic shifts $U_1$,
$U_2$ of $W_1$ and $W_2$ respectively such that $U_1 \equiv U_2$.
A subword of a cyclic word $W$ is a subword of a cyclic shift of
$W$. By $\| W \| $ we denote the length of a (cyclic) word $W$.
Finally, for a real number $r$, $[r]$ means the greatest integer
which is less than or equal to $r$.

\begin{defn}
A set $\mathcal T$ of cyclic words in $X$ satisfies the condition
$C^\ast (\lambda)$ if for all common subwords $A$ of any two
different cyclic words $B,C \in \mathcal T^{\pm 1}$, we have $\|
A\| < \lambda \min \{\| B\| , \| C\| \}$ and for all cyclic words
$B\in \mathcal T^{\pm 1}$, all subwords $A$ of $B$ of length $\|
A\| \ge \lambda \| B\| $ occur in $B$ only once.
\end{defn}

\begin{defn}
Given a set $\mathcal T$ of words in $X$, we define a {\it growth
function} of $T$ by the formula $$ \sigma _{\mathcal T} (n)=\card
\mathcal T(n),$$ where $\mathcal T(n)$ is the set of all words
from $\mathcal T$ having length exactly $n$, i.e.,
$$\mathcal  T(n)=\{ W\in \mathcal T\; :\; \| W\| =n\} .$$
\end{defn}

The main result of this section is the following.

\begin{prop}\label{WSC}
There exists a set $\mathcal T$ of words in the alphabet $X=\{a,
b\}$ and a non--increasing function $\lambda :\mathbb N \to (0,1)$
satisfying the following conditions.
\begin{enumerate}

\item[(i)] The function $\sigma_{\mathcal T}$ is non-decreasing
and $\lim_{n \to \infty}\sigma _{\mathcal T}(n)=\infty .$

\item[(ii)] $\lim_{n \to \infty} \lambda(n)= 0$.

\item[(iii)] $\mathcal T$ satisfies $C^\ast (1/50)$ condition and
for all $n\in \mathbb N$, the set $\bigcup\limits _{k=n}^\infty
\mathcal T(k)$ satisfies $C^*(\lambda (n))$.
\end{enumerate}
\end{prop}

The proof of Proposition \ref{WSC} is based on four auxiliary
lemmas. Recall that for any $l \ge 2$, a word $W$ is called
$l$-aperiodic if it has no non-empty subwords of the form $V^l$.
The following lemma can be found in the book \cite[Theorem
4.6]{Olsbook}

\begin{lem}\label{sc1}
Denote by $f(n)$ the number of all $6$-aperiodic words of length
$n>0$ over the alphabet $X=\{a, b \}$. Then we have $$ f(n)>
(3/2)^n. $$
\end{lem}

Let $\mathcal{X}(k)= \{X_{k,1},\dots, X_{k,f(k)} \}$ be the set of
all different $6$-aperiodic words of length $k$ in the alphabet
$\{a, b \}$. For every $k>8$ and every $i=0,1 \dots,
\left(\left[\frac{f(k)}{k}\right] -1\right) $, consider the
(cyclic) word

\begin{equation}
W_{k,i} = (a^6 b X_{k-8, ik+1} b)(a^6 b X_{k-8, ik+2} b) \dots
(a^6 b X_{k-8, ik+k} b). \label{Wki}
\end{equation}

 Set
$$
 \mathcal A _k= \left \{W_{k,i}:
i=0,1,\dots, \left(\left[\frac{f(k)}{k}\right] -1\right) \right\}
.
$$

The next lemma is an immediate consequence of (\ref{Wki}) and
Lemma \ref{sc1}.

\begin{lem} \label{sc2} For any $k>8$ and any $W \in \mathcal{A}_k$, we have:
\begin{enumerate}
\item[(a)] $\| W\| =k^2$;

\item[(b)] $ \# \mathcal{A}_k \ge \frac{(3/2)^k}{k} -1.$
\end{enumerate}
\end{lem}

\begin{lem}\label{SC}
For any $k>8$, the set $\bigcup \limits_{j=k}^\infty
\mathcal{A}_j$ satisfies $C^*(\frac{3}{k})$.
\end{lem}

\begin{proof}
Suppose that $U \in \mathcal{A}_j$, $j\ge k$, is a cyclic word and
$V$ is a subword of $U$ such that $\| V\| \ge (3/k) \| U\|. $ Then
we have $\| V\| \ge (3/j)\| U\| =3j>2j+8$. Note that any subword
of $U$ of length greater than $2j+8$ contains a subword of type
\begin{equation}
a^6 b X_{j-8 , i} b a^6, \label{word}
\end{equation}
where $ X_{j-8 , i} \in \mathcal X(j-8)$. Since all words from
$\mathcal {X}(j-8)$ are aperiodic and different, such a subword
occurs in $U$ only once. Therefore, $V$ occurs in $U$ once.

Further, let $U_1, U_2$ be two cyclic words from $\bigcup
\limits_{j=k}^\infty \mathcal{A}_j$ and $V$ a common subword of
$U_1, U_2$ such that $\| V\| > (3/k) \min \{\| U_1\| ,\; \| U_2\|
\} .$ Arguing as above, we can show that $V$ contains a subword of
type (\ref{word}) for $j=\min \{\| U_1\| ,\; \| U_2\| \} $. It
remains to observe that such a subword appears in a unique word
from $\mathcal{A}$.
\end{proof}

For each $n\in \mathbb N$, $n\ge 81$, we construct a set
$\mathcal{B}_n$ of words over $\{a, b \}$ as follows. First we
divide each set $\mathcal{A}_k$ into $(2k+1)$ disjoint parts such
that
\begin{equation}
\mathcal{A}_k =\bigsqcup_{i=1}^{2k+1} \mathcal{A}_{k, i},
\label{31}
\end{equation}
and
\begin{equation}
\card \mathcal{A}_{k, i} \ge \left[ \frac{\card
\mathcal{A}_{k}}{2k+1 }\right] \label{32}
\end{equation}
for any $i=1,\dots, 2k+1$. We set
\begin{equation}
\mathcal{B}_n = \mathcal{A}_{k, l} \label{33}
\end{equation}
where $k =[\sqrt{n}]$ and $l= n - [\sqrt{n}]^2$. Note that
$[\sqrt{n}]\ge [\sqrt{81}]>8$ and $ l \le n-(\sqrt{n}-1)^2 = 2
\sqrt{n} -1 \le 2k+1. $ Thus $\mathcal{B}_n$ is well-defined for
$n\ge 81$. Furthermore, for any $W\in \mathcal{B}_n$, we have
$W\in \mathcal{A}_{[\sqrt{n}]}$. Hence
\begin{equation}
n \ge |W| \ge (\sqrt{n}-1)^2 > n -2 \sqrt{n}  \label{34}
\end{equation}
by Lemma \ref{sc2}. Finally, given an arbitrary word $W\in
\mathcal{B}_n$, we form a new word
\begin{equation}
\overline{W}=W b^m \label{Wam}
\end{equation}
where $m =n -|W|$. We call $W$ a {\it core } of the word
$\overline{W}$. Inequality (\ref{34}) yields
\begin{equation}
0 \le m < 2\sqrt{n}.  \label{35}
\end{equation}
We set $$ \mathcal T(n) = \{\overline{W}: W \in \mathcal{B}_n \},
$$ for all $n\ge 81$ and $\mathcal T(n) = \emptyset$ for $n < 81$.

The proof of the next lemma is straightforward. We leave it to the
reader.

\begin{lem}\label{ComSubw}
Let $A,B,C,D$ be arbitrary words in the alphabet $X$. Suppose that
$\max \{ \| C\| , \; \| D\| \} \le y$ and any common subword of
cyclic words $A$ and $B$ has length at most $x$. Then the length
of any common subword of the cyclic words $AC$ and $BD$ is at most
$3x+2y$.
\end{lem}

\begin{proof}[Proof of Proposition 2.1.]
Let us take $\mathcal T(n)$ as defined above and set $\mathcal T=
\bigcup\limits_{k=1}^\infty \mathcal T(k)$. Combining (\ref{32}),
(\ref{33}), and Lemma \ref{sc2}, we obtain $$ \sigma(n) = \card
\mathcal T(n) = \card \mathcal{B}_n \ge \left[ \frac{\card
A_{[\sqrt{n}]}}{2[\sqrt{n}]+1} \right] \ge \left[
\frac{(3/2)^{[\sqrt{n}]}- [\sqrt{n}]} {[\sqrt{n}](2[\sqrt{n}]+1)}
\right]. $$ Evidently we have $\lim_{n \to \infty} \sigma(n) =
\infty$. Moreover, passing to a subset of $\mathcal T$ if
necessary we can always assume that $\sigma (n)$ is
non--decreasing.

Let us show that the union $\bigcup\limits_{k=n}^\infty \mathcal
T(k)$ satisfies $C^*(\lambda(n))$ for
$$ \lambda(n) =\frac{9}{[\sqrt{n}]} +
\frac{2[\sqrt{n}]}{n-2\sqrt{n}}
$$
Suppose that $\overline{U}_1$, $\overline{U}_2 $ are two different
words from $\bigcup\limits_{k=n}^\infty \mathcal T(k)$, $n\ge 81$,
and $V$ is a common subword of $\overline{U}_1, \overline{U}_2$.
Let $U_1$ and $U_2$ be the cores of $\overline{U}_1$ and
$\overline{U}_2$ respectively, $l=\min \{ \| U_1\| ,\; \| U_2\|\}
$. Note that the length of any common subword of $U_1$ and $U_2$
is at most $3l/[\sqrt{n}] $ by Lemma \ref{SC}. According to Lemma
\ref{ComSubw} and inequality (\ref{35}) this yields
$$
\frac{ \| V\| } {\min\{ \|\overline{U}_1\| , \; \|
\overline{U}_2\| \} } \le \frac{\frac{9}{[\sqrt{n}]} l
+2[\sqrt{n}]}{l}\le \lambda (n).
$$

In case $\overline{U}\in \bigcup\limits_{k=n}^\infty \mathcal
T(k)$ and $V$ is a common subword of two different cyclic shifts
of $\overline U$, we obtain the inequality $\| V\| /\| U\| \le
\lambda (n)$ in the analogous way. Finally, let $N$ be a integer
such that $\lambda (N)\le 1/50$. Then we set  $\mathcal
T(n)=\emptyset $ for all $n\le N$ and redefine $\lambda (n)$ to be
equal to $1/50$ for all $n\le N$.
\end{proof}

\section{Main construction}

Throughout the rest of the paper we fix a metric space $M$
satisfying conditions (M1) and (M2). Let $\mathcal T$ be the set
of words provided by Proposition \ref{WSC} and $\sigma =\sigma
_\mathcal T$  its growth function. Also, let us fix a sequence
$$Net_1\subseteq Net_2\subseteq \ldots $$ constructed in Lemma
\ref{net1} such that $Net_i$ is a $(2/i , 1/i)$--net in $M_i$,
$i\in \mathbb N$. By $\G _i$ we denote the complete (abstract)
graph with the vertex set $Net _i$. Further, we endow $\G_i $ with
a metric in which the length of an edge $e$ with endpoints $x$ and
$y$ is $\dm (x,y)$. Thus there is a map $\G_i\to M$ that maps each
vertex of $\G_i$ to the corresponding point of $Net _i$ and maps
edges of $\G_i$ to geodesics in $M$. It is clear that the
restriction of this map to the set of vertices of $\G_i$ is an
embedding.

We consider a sequence $n_i$, $i\in \N $, satisfying the following
three conditions:

\begin{enumerate}

\item[(I)] $\{ n_i/i\} $ is an increasing sequence of natural
numbers.

\bigskip

\item[(II)]  For any $i\in \N $, $\sigma (n_i/i)\ge N_i(N_i-1)/2,$
where $N_i=\card Net_i$.

\bigskip

\item[(III)] $n_i/i> n_{i-1}\diam M_{i-1}$ for all $i\in \N$,
$i\ge 2$.
\end{enumerate}

Obviously we can always ensure (I)--(III), choosing $n_i$ after
$n_{i-1}$. (Recall that $\sigma (n)\to \infty $ as $n\to \infty
$.)

For every $i\in \mathbb N$, we take an arbitrary orientation on
edges of $\G _i$. Let $E(\G _i)$ denote the set of all oriented
edges of $\G_i $. In the next lemma $\lceil x\rceil $ means the
smallest integer $y$ such that $y\ge x$.

\begin{lem}\label{phi}
There exists an injective labelling function $\phi
:\bigcup\limits_{i=1}^\infty E(\G_i)\to \mathcal T^{\pm 1}$ such
that for any edge $e\in E(\G_i)$ with endpoints $x,y$, we have
\begin{equation}
\| \phi (e)\| =\lceil n_i\dm (x,y)\rceil
\end{equation}
and  $\phi (e^{-1})= \phi (e)^{-1}.$ In particular, the set $\{
\phi (e)\; |\; e\in E(\G_i)\} $ of all edge labels of $\G_ i$
satisfies $C^\ast (\lambda (i))$.
\end{lem}

\begin{proof}
Suppose that $x,y\in Net _i$ and $u,v\in Net _j$ for some $j>i$.
Then combining conditions (I)--(III) and the fact that $Net_j$ is
an $(2/j , 1/j)$--net, we obtain
$$\lceil n_i\dm (x,y)\rceil\le \lceil n_i\diam M_i \rceil <\lceil
n_{i+1}/(i+1)\rceil \le \lceil n_j/j\rceil\le \lceil n_j\dm
(u,v)\rceil .$$ Thus it suffices to show that the number $l_{ik}$
of unordered pairs $x,y\in Net_i$ such that $[n_i\dm (x,y)]=k$ is
less than the number of words of length $k$ in $\mathcal T $ for
every possible $k$. Obviously we have
$$l_{ik}\le \frac{N_i(N_i-1)}{2} $$ and $$\sigma
(k)\ge \sigma (n_i/i)\ge \frac{N_i(N_i-1)}{2}$$ since $\sigma $ is
non--decreasing and $\dm (x,y)>1/i $ for any $x,y\in Net _i$.

The assertion "in particular" can be derived as follows. Note that
$n_i/i>i$ by the property (I). Since for any $e\in E(\G_i)$, we
have $$\| \phi(i)\| \ge n_i\dm (e_-, e_+)\ge n_i/i>i,$$ $\phi (e)$
belongs to the union $\bigcup\limits_{j=i}^\infty \mathcal T(j)$.
It remains to apply Proposition \ref{WSC}.
\end{proof}

If $p=e_1e_2\ldots e_n$ is a combinatorial path in $\G_i$, where
$e_1,e_2, \ldots , e_n\in E(\G_i)$, we define the label $\phi (p)$
to be the word $\phi (e_1)\phi(e_2)\ldots \phi(e_n)$. Let
$$
\mathcal R_i=\{ \phi (p)\; |\; p\; {\rm is\; an\; irreducible\;
cycle\; in\; } \G_i \}
$$
and $$\mathcal R=\bigcup\limits _{i=1}^{\infty} \mathcal R_i .$$
Finally, we define the group $G$ by the presentation
\begin{equation} \label{pres}
\langle a,b\; |\; \mathcal R\rangle .
\end{equation}

Let $\Delta $ be a van Kampen diagram over (\ref{pres}), $\Pi $ a
cell of $\Delta $. We say that $\Pi $ has {\it rank } $i$ if $Lab
(\partial \Pi )$ is a word from $\mathcal R_i$. Further, we call a
word $W$ in the alphabet $\{ a^{\pm 1}, b^{\pm 1}\} $ a {\it
$\G_i$--word} if $W$ is a label of some irreducible combinatorial
path $p$ in $\G _i$. (Evidently such a path $p$ is unique as
$\mathcal T$ satisfies $C^\ast (1/50)$ and $\phi $ is injective.)

Suppose that $p$ is a path in a van Kampen diagram $\Delta $ over
(\ref{pres}). If $Lab (p)$ is a $\G_i$--word corresponding to the
path $e_1\ldots e_t$ in $\G_i$, where $e_1, \ldots , e_t$ are
edges of $\G_i$, then $p$ can be represented as a product
\begin{equation}\label{can}
p=p_1\ldots p_t
\end{equation}
of its segments $p_1, \ldots , p_t$ with labels $Lab (p_1)=\phi
(e_1), \ldots , Lab (p_t)=\phi (e_t)$. In this case we call the
decomposition (\ref{can}) a {\it canonical} decomposition of $p$.

Note that for the boundary $p$ of a  cell $\Pi $ in $\Delta $, two edges (say, $e$ and
$f$) adjacent to the vertex $(p_i)_+=(p_{i+1})_-=e_+=f_-$ for some
$i$ can have mutually inverse labels. This allows to identify $e$
with $f^{-1}$; then we can pass to the next pair of edges adjacent
to $e_-=f_+$ and so on. Since $\mathcal T$ satisfies $C^\ast
(1/50)$, not more than $1/50 $ of each segment $p_1, \ldots , p_t
$ can be cancelled by such reductions. The irreducible path
$p_1^\prime \ldots p_t^\prime $, where $p_i^\prime $ is a subpath
of $p_i$, is called a {\it reduced boundary } of $\Pi $ and is
denoted by $\partial _{red} \Pi $. Thus we have $|p_i^\prime |\ge
\frac{48}{50} |p_i|.$ Also, to each path $q$ in $\Delta $, we
assign a path $q_{red}$ which is obtained from $q$ by eliminating
edges that do not appear in reduced boundaries of cells in $\Delta
$.

It seems more natural to consider the reduced boundary. However,
in the sequel, working with the notion of the well-attached cells
defined below, it is convenient, for technical reasons, to
distinguish between the notion of the boundary and that of the reduced
boundary.

Given two cells $\Pi _1$, $\Pi_2$ of the same rank $i$ in a van
Kampen diagram $\Delta $ over (\ref{pres}), we say that $\Pi_1$
and $\Pi_2$ are {\it well--attached to each other}, if the
following is true. Up to a cyclic shift, $\partial \Pi_1$
(respectively $(\partial\Pi_2)^{-1}$) admits a canonical
decomposition $\partial \Pi_1=p_1\ldots p_t $ (respectively
$(\partial \Pi_2)^{-1}=q_1\ldots q_s$) associated to a path
$e_1\ldots e_t$ (respectively $f_1\ldots f_s$) in $\G_i$, where
$e_1=f_1$ and $p_1=q_1$. Let $d$ be the reduced cycle in $\G_i$
obtained from $e_{2}\ldots e_tf_s^{-1}\ldots f_{2}^{-1}$. Then the
label of the cycle $c=p_{2}\ldots p_tq_s^{-1}\ldots q_{2}^{-1}$ is
freely equal to the $\G_i$--word corresponding to $d$. Thus, by
the definition of $\mathcal R_i$, $Lab(c)$ is freely equal to a
relator and hence we can replace cells $\Pi_1$ and $\Pi_2$ with
one cell (see \cite{Olsbook} for details).

Now suppose that $\partial \Delta=uw$, where $Lab (w)$ is a $\G
_i$--word. We say that a cell $\Pi $ of rank $i$ is {\it
well--attached to a segment $w$ of boundary of $\Delta $} if, up
to a cyclic shift, $\partial \Pi $ (respectively $w^{-1}$) admits
a canonical decomposition $\partial \Pi =p_1\ldots p_t $
(respectively $w^{-1}=q_1\ldots q_s$) associated to a path
$e_1\ldots e_t$ (respectively $f_1\ldots f_s$) in $\G_i$, where
$e_1=f_1$ and $p_1=q_1$. In this case we denote by $d$ the reduced
cycle in $\G_i$ obtained from $f_s^{-1}\ldots f_{2}^{-1}e_2\ldots
e_t$. Then the label of the path $v=q_s^{-1}\ldots
q_{2}^{-1}p_{2}\ldots p_t$ is freely equal to the $\G_i$--word
corresponding to $d$. Thus, by cutting the cell $\Pi $, we obtain
a subdiagram $\Sigma $ of $\Delta $ such that $\partial \Sigma
=uv$, where $v$ is also a $\G_i$--word.

We can summarize these observations as follows.

\begin{lem}\label{wa1}
Let $\Delta $ be a van Kampen diagram over (\ref{pres}).
\begin{enumerate}
\item Suppose that $\Delta $ is minimal, i.e., it has minimal
number of cells among all diagrams over (\ref{pres}) with the same
boundary label. Then no two cells of $\Delta $ are well attached
to each other.

\medskip

\item Suppose that $Lab (\partial \Delta )=VW$, where $W$ is a
$\G_i$--word. Assume that a cell $\Pi $ is well--attached to the
subpath of $\partial \Delta $ labelled $W$. Then there exists a
subdiagram $\Sigma $ of $\Delta $ (which can be obtained from
$\Delta $ by cutting the cell $\Pi $) such that $Lab (\partial
\Sigma) =VU$, where $U$ is a $\G_i$--word.
\end{enumerate}
\end{lem}

The next lemma provides certain sufficient conditions for two
cells (or a cell and a part of boundary of a diagram) to be
well--attached.

\begin{lem}\label{wa2}  Let $\Delta $ be a van Kampen diagram over
(\ref{pres}).

\begin{enumerate}

\item Suppose that $\Pi _1$, $\Pi_2$ are cells in $\Delta $ such
that there exists a common subpath $p$ of $\partial _{red}\Pi_1$
and $(\partial _{red}\Pi_2)^{-1}$ such that
$$
|p|\ge \frac{1}{10} \min\{ |\partial _{red}\Pi_1 |,\; |\partial
_{red} \Pi_2| \}.
$$
Then $\Pi_1$ and $\Pi_2$ are well--attached to each other.

\medskip

\item Suppose that $\partial \Delta =vw$, where $Lab (w)$ is a
$\G_i$--word. Assume that for a cell $\Pi $, there is a common
subpath $q$ of $w_{red}$ and $(\partial _{red} \Pi )^{-1}$ of
length $$ |q|\ge \frac{1}{10} |\partial _{red}\Pi |.$$ Then $\Pi $
is well--attached to the subpath $w$ of $\partial \Delta $.
\end{enumerate}
\end{lem}

\begin{proof}
Let us prove the first assertion of the lemma. Up to a cyclic
shift, $\partial \Pi_1$ admits canonical decomposition $p_1\ldots
p_t$. Let $p_1^\prime \ldots p_t^\prime $ be the corresponding
decomposition of $\partial _{red}\Pi_1 $, where $p_i^\prime $ is a
subpath of $p_i$. Then $p$ and a certain $p_i^\prime $ have a
common subpath $q$ of length at least $(1/20) |p_i^\prime|\ge
(48/1000) |p_i|$. Let $q_1\ldots q_s$ be the canonical
decomposition of $\partial \Pi_2$, $q_1^\prime \ldots q_s^\prime $
the corresponding decomposition of $\partial _{red}\Pi_2$. Let $Z$
denote the set of endpoints of paths $q_1^\prime , \ldots ,
q_s^\prime $. If $q$ is cut by vertices from $Z$ into at most two
parts, then one of these parts has length at least $1/2 |q|\ge
(24/1000) |p_i|>(1/50) |p_i|$. If $q$ contains more than one
vertex from $Z$, then $q_j^\prime $ is a subpath of $p_i$ for some
$j$ (note that $|q_j^\prime |\ge (48/50) |q_j|$). In both cases we
found a common subpath of $p_i$ and $q_j$ of length at least
$(1/50) \min\{ |p_i|,\; |q_j|\} $. Therefore the labels of edges
$e_i$ and $f_j$ of $\G_k $ and $\G_l$ respectively corresponding
to $p_i$ and $q_j$ contain a common subword of length at least
$(1/50) \min\{ \|\phi (e_i)\| ,\; \|\phi (f_j)\| \} $. Since
$\mathcal T$ satisfies $C^\ast (1/50)$ and $\phi $ is injective,
we have $p_i=q_j$, $k=l$, and $e_i=f_j$. The proof of the second
assertion is similar and we leave it to the reader.
\end{proof}

From Lemmas \ref{wa1} and \ref{wa2}, we immediately obtain

\begin{cor}\label{scc}
\begin{enumerate}
\item  Let $\Delta $ be a minimal van Kampen diagram over
(\ref{pres}). Then for any common subpath $p$ of the reduced
boundaries any two cells $\Pi _1$ and $\Pi _2$ of $\Delta $, we
have $|p|< \frac{1}{10} \min\{ |\partial _{red}\Pi_1 |,\;
|\partial _{red}\Pi_2| \}. $
\end{enumerate}
\end{cor}

Up to notation, the proof of the next lemma coincides with the
proof of Lemma 8 in \cite{Ols}. We provide it for convenience of
the reader.

\begin{lem}\label{VW}
Suppose that $W$ is a $\G_i$--word. Then there exists a word $V$
such that $W=V$ in $G$,  $V$ is of the minimal length among all of
the words (not necessarily $\G_i$--words) representing the same
element as $W$ in $G$, and $V$ is freely equal to a $\G_i$--word.
\end{lem}

\begin{proof}
Let $V$ be a shortest word representing the same element as $W$ in
$G$. We consider a van Kampen diagram $\Delta $ over (\ref{pres})
corresponding to this equality. Without loss of generality we may
assume that the word $W$ and $\Delta $ are chosen in such a way
that $\Delta $ has the minimal number of cells among all diagrams
corresponding to equalities of $V$ to $\G_i$--words. We are going
to show that $\Delta $ contains no cells at all, and thus $V$ is
freely equal to a $\G_i $--word.

Assume that there is at least one cell in $\Delta $. Denote by
$\Delta ^\prime $ the map obtained from $\Delta $ by eliminating
all edges that do not appear in reduced boundaries of cells of
$\Delta $. Then $\Delta ^\prime $, as a map, satisfies $C^\prime
(1/10)$ small cancellation condition (see \cite{LS}[Chapter 5]) by
Corollary \ref{scc}. By Greendlinger's Lemma, this means that
$\Delta $ contains a cell $\Pi $ such that there is a common
subpath of $\partial _{red}\Pi $ and $(\partial\Delta )_{red}$ of
length $|p|>0.7 |\partial _{red}\Pi |.$ (We substitute $\lambda
=0.1$ in the Greendlinger's constant $1-3\lambda $ from
\cite{LS}).

The boundary of $\Delta $ consists of two parts $v$ and $w$
corresponding to words $V$ and $W$. If the path $p$ has a common
subpath with $w^{\pm 1}_{red}$ of length at least $0.1|\partial
\Pi |$, then $\Pi $ is well--attached to the subpath $w$ of
$\partial \Delta $ by Lemma \ref{wa2}. However, by the second
assertion of Lemma \ref{wa1} this contradicts to the choice of $W$
and $\Delta $. Hence there is a common subpath $q$ of $\partial
_{red} \Pi $ and $v^{\pm 1}_{red}$ such that $|q|>
(0.7-0.2)|\partial _{red}\Pi |=0.5|\partial _{red}\Pi |$. Thus
$\partial _{red}\Pi =qq_1$, $|q|>|q_1|$, and the words $Lab (q)$,
$Lab (q_1)$ represent the same element in the group $G$. But $Lab
(q)$ is a subword of $V$ and we arrive at a contradiction to our
choice of $V$ as a shortest word representing the same element as
$W$ in $G$.
\end{proof}

\begin{defn}
For each $i\in\mathbb N$, we construct an embedding $$\alpha _i:
Net_i \to G $$ as follows. Let us fix a point $O$ in $Net_1$ (and
thus $O\in Net_i$ for all $i$). Then for any $x\in Net_i$ there is
a combinatorial path $p$ in $\G_i$ such that $p_-=O$, $p_+=x$. We
define $\alpha _i(x)$ to be equal to the element of $G$
represented by $\phi (p)$. Note that $\alpha _i(x)$ is independent
of the choice of $p$. Indeed, if $q$ is another path in $\G_i$
with the origin $O$ and terminus $x$, then $pq^{-1}$ is a cycle
and thus $\phi(p)\phi(q^{-1})$ is a relator from $\mathcal R_i$,
i.e., $\phi (p)$ and $\phi (q)$ represent the same element of $G$.
\end{defn}

\begin{lem}\label{emb1}
Let $\dg $ denote the word metric on $G$ corresponding to the
generating set $\{ a,b\}$. Then for any $i\in\mathbb N$ and any
$x,y\in Net_i$, we have
\begin{equation}\label{dgdm}
(1-2\lambda(i)) \dm (x,y)\le \frac{1}{n_i} \dg (\alpha
_i(x),\alpha _i(y))\le \dm (x,y)+ \frac{1}{n_i}.
\end{equation}
\end{lem}

\begin{proof}
If $e$ is an edge in $\G_i$  such that $e_-=x$, $e_+=y$, and $a,b$
are edges in $\G_i $ such that $a_-=b_+=O$, $a_+=x$, $b_-=y$, then
$aeb$ is a cycle in $\G_i$. Therefore $\phi(a)\phi(e)\phi(b)$
labels a cycle $c$ in $\Ca $ with beginning at $1$. Let $c=psq$,
where $Lab (p)\equiv \phi (a)$, $Lab (s)\equiv \phi (e)$, $Lab
(q)\equiv \phi (b)$. Since by definition $p_+=\alpha _i(x)$ and
$q_-=\alpha _i(y)$, the elements $\alpha _i(x)$ and $\alpha _i(y)$
are connected by the path $s$ in $\Ca $. Therefore,
$$
\dg (\alpha_i(x),\alpha_i(y))\le |s|=\| \phi (e)\| =\lceil n_i \dm
(x,y)\rceil \le n_i\dm (x,y) +1.
$$
This gives the right hand side inequality in (\ref{dgdm}).

Further, by Lemma \ref{VW}, there exists a word $V$ representing
the element $(\alpha_i(x))^{-1}\alpha_i(y))$ and a $\G_i$--word
$U$ freely equal to $V$ such that
\begin{equation} \label{dgdm0}
\| V\| = \dg (1, \alpha_i(x))^{-1}\alpha_i(y))= \dg (\alpha_i(x),
\alpha_i(y)).
\end{equation}
Obviously we have
\begin{equation}\label{dgdm1}
\| V\| \ge (1-2\lambda (i))\| U\|
\end{equation}
since the set of edge labels of $\G_{i}$ satisfies $C^\ast
(\lambda (i))$ by Lemma \ref{phi}. Let $r=e_1\ldots e_t$ be the
path in $\G_i$ corresponding to $U$. Then, arguing as in the first
case, we can show that $r_-=x$, $r_+=y$ and thus
\begin{equation}\label{dgdm2}
\|U\| = \sum\limits_{j=1}^t \| \phi (e_j)\| =\sum\limits_{j=1}^t
n_i\dm ((e_j)_-,(e_j)_+)\ge n_i\dm (x,y).
\end{equation}
Combining (\ref{dgdm0}), (\ref{dgdm1}), and (\ref{dgdm2}) we
obtain the left hand inequality in (\ref{dgdm}).
\end{proof}

\begin{defn}
We take a non--principal ultrafilter $\omega $ such that $\omega
(\{ n_i\} )=1$ and consider the asymptotic cone $\C $ of $G$ with
respect to this ultrafilter. Our next goal is to define an
embedding $\iota $ of $M$ to $\C $.

Let $x$ be a point of $M$. Then there is a sequence of points
$x_i\in Net_i$ such that $x_i\to x$ as $i\to \infty $. We define
$\iota (x)$ to be the point of $\C $ represented by an arbitrary
sequence $\{ g_i\} $, where  $g_{n_i}=\alpha_i (x_i)$ for any
$i\in \mathbb N$. Obviously $\iota $ is well--defined as the point
of $\C $ representing the sequence $\{ g_i\} $ depends on the
subsequence $\{ g_{n_i}\} $ only.
\end{defn}

\begin{prop} Suppose that $M$ is a metric space satisfying (M1)
and (M2). Then the map $\iota $ is an isometry.
\end{prop}

\begin{proof}
Let $x,y$ be points of $M$, $\{ x_i \}$, $\{ y_i\} $ the sequences
of elements of nets $Net_i$ such that $x_i\to x$ as $i\to \infty $
and $y_i\to y$ as $i\to \infty $. Let $\{ g_i\} $ and $\{ h_i\} $
be the corresponding sequences of elements of $G$ representing
$\iota (x)$ and $\iota (y)$. Then applying Lemma \ref{emb1}, we
have
$$
\begin{array}{rl}
\dm (x,y)= & \lim\limits_{i\to\infty }
\dm(x_i,y_i)=\lim\limits_{i\to\infty } \frac{1}{n_i} \dg (\alpha
_i(x_i),\alpha _i(y_i))= \\ & \\ & \lim\limits_{\omega }
\frac{1}{i} \dg (g_i,h_i)=\dc (\iota (x), \iota (y)).
\end{array}
$$
\end{proof}

\section{Embedding of the fundamental group}
All assumptions and notation from the previous section remain in
force here. In particular, $\iota $ denotes the isometry $M\to \C
$ constructed in the previous section. In addition we suppose that
$M$ satisfies (M3). Also, let $\iota ^\ast $ denote the
homomorphism $\pi_1(M)\to \pi_1(\C )$ induced by $\iota $. We
conclude the proof of Theorem \ref{main} by proving the following.

\begin{prop}
Suppose that $M$ satisfies (M1)--(M3). Then the map $\iota ^ \ast
:\pi_1(M)\to \pi_1(\C )$ is injective.
\end{prop}

\begin{proof}

Let $S= [0,1] \times [0,1]$ be a unit square and $\gamma:
\partial S \to M$ a loop in $M$ such that $\iota\gamma $ is
contractible in $\C $. We want to show that $\gamma $ is
contractible in $M$.

Since $\iota\gamma $ is contractible in $\C $, there exists a
continuous map  $r: S \to \C $ such that the restriction of $r$ to
$\partial S$ coincides with $\iota\gamma $. The unit square $S$ is
compact, and therefore $r$ is uniformly continuous. Hence there
exists $\delta $ such that for any $y_1, y_2 \in B$ which lie at
distance at most $\delta $ in $B$, we have
\begin{equation}\label{40}
\dc(r(x), r(y)) < \e/20.
\end{equation}
We can also assume that $1/\delta \in \N$. By $\gr $ we denote the
standard $\delta $--net in $S$ that is the set $$\gr =\{ (a
\delta, b \delta)\; |\; a, b \in \Z,\; 0 \le a, b \le 1/\delta \}
.$$ By $r (\gr )$ we denote the image of $\gr $ in $\C $.

For every point $x\in r(\gr )\cup\left(
\bigcup\limits_{i=1}^\infty \iota (Net_i) \right) $ we fix an
arbitrary sequence $\{ x_i\} $ of elements of $G$ that represents
$x$ in $\C $ (such a sequence will be called a {\it standard
representative} of $x$). Let $\e $ be the constant from (M3).  We
take $L\in \mathbb N$ such that the following conditions hold:
\begin{enumerate}
\item[(L0)] $r(\gr )$ is contained in $\iota (M_L)$;

\item[(L1)] $ 1/L< \e /20 $; in particular, $1/n_L< \e/20 $;

\medskip

\item[(L2)] for any two points $x,y\in r(\gr )\cup \iota (Net_L)$,
we have
$$ \left| \frac{1}{n_L} \dg (x_{n_L}, y_{n_L}) - \dc (x,
y)\right| \le \e /20 , $$ where $\{ x_i\}$, $\{ y_i\} $ are
standard representatives of $x$ and $y$ respectively.
\end{enumerate}

(Note that for any $l$ there exist $L>l$ such that  (L0)--(L2)
hold.)

We say that two points $x,y$ in $\gr $ are neighbors, if they have
the form $x=( a \delta, b \delta)$, $y=( (a+1)\delta, b \delta)$
or $x=( a \delta, b \delta)$, $y=( a\delta, (b+1) \delta)$. If
$x,y\in \gr$ are neighbors and $\{ x_i\} $ , $\{ y_i\} $ are
standard representatives of $r(x),r(y)$, we fix an arbitrary
geodesic in the Cayley graph $G$ going from the element $x_{n_L}$
to $y_{n_L}$ and denote this geodesic by $g(x_{n_L},y_{n_L})$.
Further for every point $x\in \gr $ which lies on $\partial S $,
we take a point $t^x\in \iota (Net _L)$ which is closest to
$r(x)$; in particular, we have
\begin{equation} \label{50}
\dc (t^x, r(x))\le 2/L \le 0.1\e
\end{equation}
as $r(x)\in \iota (M_L)$ by (L0) and $\iota (Net_L) $ is a
$2/L$--net in $\iota (M_L)$ (recall that $\iota $ is an isometry).
Suppose that $\{ t^x_i\} $ is the standard representative of
$t^x$. Then we join elements $x_{n_L}$ and $t^x_{n_L}$ by a
geodesic $h(x_{n_L}, t^x_{n_L})$ in $\Ca $. Finally, if $x,y\in
\partial S$ are neighbors and $t^x,t^y$ are the corresponding
points of $\iota (Net_L)$, then we denote by $k(t^x_{n_L},
t^y_{n_L})$ a path in $\Ca $ joining $t^x_{n_L}$ to $t^y_{n_L}$
such that the label of $k$ is equal to $\phi (e)$, where $e$ is
the edge of $\G_L$ satisfying the conditions $\iota (e_-)=t^x$,
$\iota (e_+)=t^y$. In particular, we have
\begin{equation} \label{55}
|k(t^x_{n_L}, t^y_{n_L})| =\| \phi (e)\| =\lceil n_L\dc (t^x, t^y)
\rceil
\end{equation}

Let $x^1, \ldots , x^{m}$, where $m=4/\delta $, be subsequent
points of $\gr \cap\partial S$ (i.e., $x^i$ and $x^{i+1}$ are
neighbors, where indices are modulo $m$). Then the label of the
cycle
\begin{equation}\label{boundary}
p=k(t^{x^1}_{n_L}, t^{x^2}_{n_L})k(t^{x^2}_{n_L},
t^{x^3}_{n_L})\ldots k(t^{x^m}_{n_L}, t^{x^1}_{n_L})
\end{equation}
is a $\G_i$--word. We construct a van Kampen diagram $\Xi $ with
boundary label $Lab (\partial \Xi )\equiv \phi (p)$ as follows.
The net $\gr $ allows to regard $S$ as a union of $\frac{1}{\delta
^2}$ small squares with sides of length $\delta $. For any such a
square with vertices $x, y, z, t$ in $\gr $, we consider a minimal
van Kampen diagram (homeomorphic to a disk) with boundary label
\begin{equation}\label{type1}
Lab(g(x_{n_L},y_{n_L})g(y_{n_L},z_{n_L})g(z_{n_L},t_{n_L})g(t_{n_L},x_{n_L})).
\end{equation}
Also, if $x,y\in \partial S\cap \gr $ are neighbors, we consider a
minimal van Kampen diagram (homeomorphic to a disk) with boundary
label
\begin{equation}\label{type2}
Lab(g(x_{n_L},y_{n_L})h(y_{n_L},t^y_{n_L})k(t^y_{n_L},t^x_{n_L})
(h(x_{n_L},t^x_{n_L}))^{-1}).
\end{equation}
We call the constructed diagrams with boundary labels
(\ref{type1}), (\ref{type2}) {\it elementary}. Gluing these
elementary diagrams together in the obvious way we obtain a
diagram $\Xi $ over (\ref{pres}) such that $Lab (\partial \Xi )$
is the $\G _i$--word defined by (\ref{boundary}).

We are going to show that the perimeter of each elementary diagram
is less than $0.7n_L\e $. Indeed, inequality (\ref{40}) and
condition (L2) together yield
\begin{equation}\label{60}
|g(x_{n_L}, y_{n_L})|=\dg (x_{n_L}, y_{n_L})\le n_L\left( \dc
(r(x), r(y))+0.05\e\right) \le 0.1n_L\e
\end{equation}
for any two neighbors $x,y\in \gr $. If $x\in \gr \cap \partial S
$, then (L1),(L2) and (\ref{50}) imply
\begin{equation}\label{70}
\begin{array}{l}
|h(x_{n_L},t^x_{n_L})|=\dg (x_{n_L}, t^x_{n_L}) \le n_L\left( \dc
(r(x), t^x)+0.05\e\right) \\ \\ \le n_L (2/L +0.05\e )\le
0.15n_L\e .
\end{array}
\end{equation}
Finally, if $x,y\in \partial S$ are neighbors, then combining
(\ref{40}), (\ref{50}), and (\ref{55}) we obtain
\begin{equation}\label{70}
\begin{array}{l}
|k(t^x_{n_L},t^y_{n_L})|=  \lceil n_L\dc (t^x,t^y)\rceil \le \\  \\
 \lceil n_L\left(\dc (t^x, r(x))+\dc (r(x), r(y))+\dc
(r(y),t^y)\right)\rceil \le \\  \\  \lceil 0.25n_L\e \rceil \le
0.25n_L\e +1\le 0.3n_L\e
\end{array}
\end{equation}
Therefore, any word of type (\ref{type1}) or (\ref{type2}) has
length at most $0.7 n_L\e $.

\begin{lem}
Let $\Pi$ be a cell of rank $L$ in $\Xi $, $l$ the loop in $\G_L$
corresponding to the $\G_L$--word $Lab (\partial \Pi )$. Then $l$
is contractible in $M$.
\end{lem}

\begin{proof} Note that $\Pi$ lies in some elementary diagram $\Theta $.
Since any elementary diagram is minimal, it satisfies $C'(1/10)$
small cancellation condition as a map by Corollary \ref{scc}.
Hence the length of the reduced boundary of any cell in $\Theta $
is not greater than $|\partial \Theta |\le 0.7n_L\e $. This means
that
$$ |\partial \Pi| \le \frac{50}{48}| \partial _{red}\Pi | \le
\frac{35}{48}n_L\e <n_L\e .$$ Let $l=e_1\ldots e_t$, where $e_1,
\ldots , e_t$ are edges of $\G_L$. The length of $l$ satisfies
$$ |l|= \sum\limits_{i=1}^t |e_i|\le \sum\limits_{i=1}^t \frac{1}{n_L}\|
\phi (e_i)\| =\frac{1}{n_L}|\partial \Pi | < \e .$$ Therefore, $l$
is contractible in $M$ by (M3).
\end{proof}

\begin{lem}\label{q}
Consider a van Kampen diagram $\Delta $ with boundary labelled by
a $G_L$--word. Suppose that boundary label of each  cell of rank
$L$ in this diagram corresponds to a contractible loop in $M$.
Then the boundary label of the diagram also corresponds to a
contractible loop in $M$.
\end{lem}

\begin{proof}
We prove the statement of the lemma by induction on the number $s$
of cells in the diagram. If $s=0$ the statement is obvious, so we
assume that $s\ge 1$.

By Grindlinger's lemma at least one of the following two
statements holds.

1) There exist two cells $\Pi_1$ and $\Pi_2$  and a common subpath
$p$ of $\partial _{red}\Pi_1 $ and $\partial _{red}\Pi_2 $ such
that $|p|\ge \frac{1}{10} \min\{ |\partial _{red}\Pi_1 |,\;
|\partial _{red}\Pi_2| \} $.

2) There exist a cell $\Pi $ and a common subpath $p$ of $\partial
_{red} \Pi $ and $\partial \Delta $ such that $|p|\ge \frac{7}{10}
|\partial _{red}\Pi | $.

In the first case $\Pi _1$ and $\Pi _2$ have the same rank and are
well--attached to each other by Lemma \ref{wa2}. Arguing as in the
proof of Lemma \ref{wa1}, we can replace $\Pi_1 $ and $\Pi_2$ by
one cell $\Upsilon $. If $rank\; \Pi_1=  rank\; \Pi_2\ne L$, the
statement is true by the inductive hypothesis. To use the
inductive hypothesis in case $rank\; \Pi_1= rank\; \Pi_2=L$, we
have to check that the cycle $s$ corresponding to the new cell
$\Upsilon $ is contractible in $M$. Indeed, if $p,q$ are cycles
corresponding to $\Pi_1$ and $\Pi_2$ (we may assume that
$p_-=q_-$), then $s$ is homotopic to the product of $p$ and
$q^{-1}$. Since $p$ and $q$ are contractible in $M$ by the
condition of the lemma, $s$ is contractible in $M$.

In the second case $\Pi $ has rank $L$ by Lemma \ref{wa1} and is
well--attached to the boundary of $\Delta $. We pass to the
subdiagram $\Sigma $ of $\Delta $ obtained by cutting the cell
$\Pi $. Applying Lemma \ref{wa1} again, we conclude that $Lab
(\partial \Sigma )$ is a $\G_L$--word. By the inductive assumption
the cycle $c$ corresponding to $Lab(\partial \Sigma )$ is
contractible in $M$. Let $d$ be the cycle in $\G_L$ corresponding
to $Lab (\partial \Pi )$, $f$ the cycle corresponding to $Lab
(\partial \Delta )$ . As in the previous case, $f$ is homotopic to
the product of $c$ and $d^{-1}$ and hence is contractible in $M$.
\end{proof}

Now we return to the proof of the proposition. The two previous
lemmas imply that the loop $q$ in $\G_L$, corresponding to the
boundary label of the diagram $\Xi $ under consideration is
contractible.

As above, let $x_1, \ldots , x_m$ be subsequent neighbors in $\gr
\cap \partial S$. For every two neighbors $x_i,x_{i+1}$ (indices
are modulo $m$), we denote by $c_i, d_i, e_i, f_i$ the segment
$[r(x_i), r(x_{i+1})]$ of $\iota \gamma $, the geodesic path from
$r(x_{i+1})$ to $t^{x_{i+1}}$, the edge $e$ of $\G _i$ such that
$\iota (e_-)=t^{x_{i+1}}$, $\iota (e_+)= t^{x_{i}}$, and the
geodesic path from $t^{x_{i}}$ to $r(x_{i})$ respectively. Note
that for any $i=1, \ldots , m$, the the cycle $a_i=c_id_ie_if_i$
is contained in the ball $B_i=B_i(0.35\e , x_i)$ of radius $0.35\e
$ around $r(x_i)$ in $\C $. Indeed any point of $c_i$ is contained
in $B_i$ by (\ref{40}). Further since $d_i$ and $f_i$ are
geodesic, $d_i$ and $f_i$ are contained in $0.1\e $--neighborhoods
of $r(x_{i+1})$ and $r(x_i)$ respectively according to (\ref{50});
together with (\ref{40}) this implies that $d_i$ and $c_i$ lay in
$B_i$. Finally, each point of $e_i$ belongs to $B_i$ as $e_i$ is
geodesic, the distance between $r(x_i)$ and the end of $e_i$ is at
most $0.1\e $, and the length of $e_i$ is at most $0.25\e $ by the
triangle inequality. Thus $a_i$ is contained in $B_i$. Since
$\iota $ is an isometry, this means that the preimage of $a_i$
under $\iota :M\to \C$ is contractible in $M$ by (M3). Hence
$\iota (q)$ is homotopic to $\iota\gamma $ via a homotopy in $M$.
Hence $\gamma$ is contractible in $M$ according to Lemma \ref{q}.
\end{proof}

\section{Concluding remarks and questions}

We have shown that any countable group can be embedded into a
fundamental group of an asymptotic cone of some finitely generated
group. Note that our proof also shows that any recursively
presentable group can be embedded into a fundamental group of some
finitely presentable group.

The construction of our group depends on a  space $M$ and a
scaling sequence $n_k$. Similarly we can start with  a countable
set of  spaces $N_j$ (satisfying (M1) -(M3)), take a countable set
of non-intersecting scaling sequences $n_k^j(N_j)$ and construct a
group $G$, such that for each $j$ there is a scale on which $N_j$
is embedded into the asymptotic cone of $G$. A natural task is to
check that starting with the spaces with very different
fundamental groups (e.g. $\Z/ p\Z$ for different $p$) one gets
asymptotic cones (on different scales) with infinitely many
different fundamental groups. Then under certain conditions on the
spaces  the group $G$ is recursively presentable and we can embed
it into a finitely presentable group. Again, a natural task is to
check that one can chose this embedding in such a way that this
finitely presented group has different fundamental groups on
different scales.

Another natural question is: does there exists a finitely
presented group such that the simple connectivity of the
asymptotic cone depends on the choice of the ultrafilter?

Finally let us mention that  recently L.Kramer,S. Shelah, K.Tent
and S. Thomas \cite{new} have shown that if continuum hypothesis
fails, than there exist finitely presented groups (which are
uniform lattices in certain semisimple Lie groups) that have
infinitely many different asymptotic cones. However, if continuum
hypothesis holds, than the examples from \cite{new} have unique
asymptotic cones.

{\bf Acknowledgements.} The results of this paper were obtained in
the Spring 2001, as both authors were visiting the University of
Geneva. We would like to thank Pierre de la Harpe for his
invitation to Geneva and to the Swiss National Science Foundation
for the support of our work. We are grateful to the referee for
useful remarks and suggestions. We are also grateful to Mark Sapir and Cornelia Drutu
for the interest they have shown to our results and 
 Alexander Ol'shanskii for helpful comments and remarks.

\end{document}